\documentclass[a4]{article}

\usepackage{amsmath}
\usepackage{amssymb}
\usepackage{amsthm}
\usepackage{amsxtra}
\usepackage{pstcol}
\usepackage{graphicx}
\usepackage{ifthen}

\def\de{Definition}
\def\pp{Proposition}
\def\th{Theorem}

\def\lm{Lemma}
\def\co{Corollary}
\def\ex{Example}
\def\exs{Examples}
\def\re{Remark}

\newcommand{\comment}[1]{}

\newcommand{\is}{\ensuremath{\cap}}

\renewcommand{\vert}{\ensuremath{\,|\,}}
\newcommand{\norm}[1]{\ensuremath{\vert\!\!\vert{#1}\vert\!\!\vert}}

\renewcommand{\deg}[2][{}]{\ifthenelse{\equal{#1}{}}
				{\normalfont\ensuremath{\textrm{deg}(#2)}}
				{\normalfont\ensuremath{\textrm{deg}_{\mathcal #1}(#2)}}}

\newcommand{\ordP}[1]{\normalfont\ensuremath{\textrm{ord$_P$}(#1)}}
\newcommand{\Spec}[1]{\normalfont\ensuremath{\textrm{Spec}(#1)}}
\newcommand{\Max}[1]{\normalfont\ensuremath{\textrm{Max}(#1)}}
\newcommand{\Hess}[1]{\normalfont\ensuremath{\textrm{Hess}_#1}}
\newcommand{\rank}[1]{\normalfont\ensuremath{\textrm{rank}(#1)}}
\newcommand{\Sing}[1]{\normalfont\ensuremath{\textrm{Sing}(#1)}}

\renewcommand{\det}[1]{\normalfont\ensuremath{\textrm{det}\!\left(#1\right)}}
\newcommand{\Rad}[1]{\normalfont\ensuremath{\textrm{Rad}(#1)}}

\newcommand{\nfrac}[2]
{\raisebox{0.5ex}{\footnotesize \ensuremath{#1}}\hspace{-0.1ex}\raisebox{0.1ex}{/}\hspace{-0.2ex}\raisebox{-0.4ex}{\footnotesize \ensuremath{#2}}}

\newcommand{\ifdef}[2]{\ifthenelse{\boolean{#1}}{#2}{}}

\def\hA{\section}

\def\hAA{\section*}

\def\Index#1{}

\newcommand{\headA}[1]{\hA{#1}\Index{#1}}

\newcommand{\headAA}[1]{\hAA{#1}\Index{#1}\addcontentsline{toc}{section}{#1}
}

\theoremstyle{definition}
\newtheorem{Definition}{\de}[section]
\newtheorem{Example}[Definition]{\ex}
\newtheorem{Examples}[Definition]{\exs}
\newtheorem{Remark}[Definition]{\re}

\theoremstyle{plain}
\newtheorem{Proposition}[Definition]{\pp}
\newtheorem{Lemma}[Definition]{\lm}
\newtheorem{Corollary}[Definition]{\co}
\newtheorem{Theorem}[Definition]{\th}

\newcommand{\defn}[2][{}]{\ifthenelse{\equal{#1}{}}
				{\begin{Definition}#2\end{Definition}}
				{\begin{Definition}[#1]#2\end{Definition}}}
\newcommand{\prop}[2][{}]{\ifthenelse{\equal{#1}{}}
				{\begin{Proposition}#2\end{Proposition}}
				{\begin{Proposition}[#1]#2\end{Proposition}}}
\newcommand{\lma}[2][{}]{\ifthenelse{\equal{#1}{}}
				{\begin{Lemma}#2\end{Lemma}}
				{\begin{Lemma}[#1]#2\end{Lemma}}}
\newcommand{\cor}[2][{}]{\ifthenelse{\equal{#1}{}}
				{\begin{Corollary}#2\end{Corollary}}
				{\begin{Corollary}[#1]#2\end{Corollary}}}
\newcommand{\thm}[2][{}]{\ifthenelse{\equal{#1}{}}
				{\begin{Theorem}#2\end{Theorem}}
				{\begin{Theorem}[#1]#2\end{Theorem}}}
\newcommand{\rem}[2][{}]{\ifthenelse{\equal{#1}{}}
				{\begin{Remark}#2\end{Remark}}
				{\begin{Remark}[#1]#2\end{Remark}}}
\newcommand{\expl}[2][{}]{\ifthenelse{\equal{#1}{}}
				{\begin{Example}#2\end{Example}}
				{\begin{Example}[#1]#2\end{Example}}}
\newcommand{\expls}[2][{}]{\ifthenelse{\equal{#1}{}}
				{\begin{Examples}#2\end{Examples}}
				{\begin{Examples}[#1]#2\end{Examples}}}
\newcommand{\prf}[2][{}]{\ifthenelse{\equal{#1}{}}
				{\begin{proof}#2\end{proof}}
				{\begin{proof}[#1]#2\end{proof}}}

\begin{document}

\begin{center}
{\huge\bf The space of nodal curves of type $p,q$ with given  Weierstra\ss\ semigroup}
\end{center}
\begin{center}
	{H.\ Knebl}
	\footnote{Fakult\"at Informatik, Georg-Simon-Ohm-Hochschule N\"urnberg, Kesslerplatz 12, 90489 N\"urnberg, Germany},
	{E.\ Kunz $^2$ and R.\ Waldi}
	\footnote{Fakult\"at f\"ur Mathematik, Universit\"at Regensburg, Universit\"atsstrasse 31, 93053 Regensburg, Germany}
\end{center}
\vspace{0.5cm}

\begin{abstract}
We continue the investigation of curves of type $p,q$ started in [KKW]. We study the space of such curves and the space 
of nodal curves with prescribed Weierstra\ss\ semigroup. A necessary and sufficient criterion for a numerical semigroup 
to be a Weierstra\ss\ semigroup is given. 
We find a class of Weierstra\ss\ semigroups which apparently has not yet been described in the literature.
\end{abstract}

\headAA{Introduction}

Let $K$ be an algebraically closed field of characteristic $0$. For relatively prime numbers $p,q\in \mathbb N$ with 
$1<p<q$ a plane curve $C$ of type $p,q$ is the zero-set of a Weierstra\ss\ polynomial of type $p,q$
$$F(X,Y):=Y^p+bX^q+\sum_{\nu p+\mu q<pq}b_{\nu\mu}X^{\nu}Y^{\mu}\ (b_{\nu\mu}\in K, b\in K\setminus\{0\})$$
in $\mathbb A^2(K)$.
Such curves are irreducible and have only one place $P$ at infinity, i.e.\ $P$ is the only point at infinity of 
the normalization of the projective closure $\mathcal R$ of $C$. The Weierstra\ss\ semigroup $H(P)$ of $\mathcal R$ 
at $P$ is also called the Weierstra\ss\ semigroup of $C$. It contains the semigroup $H_{pq}$ generated by $p$ and $q$ 
as a subsemigroup. Hence $H(P)$ is obtained from $H_{pq}$ by closing some of its $d:= \nfrac{1}{2}(p-1)(q-1)$ gaps. 
Remember that $H_{pq}$ is a symmetric semigroup with conductor $c:=(p-1)(q-1)$. It is shown in [KKW] that any 
Weierstra\ss\ semigroup is the Weierstra\ss\ semigroup of a plane curve of type $p,q$ having only nodes as 
singularities if $p$ and $q$ are properly chosen.

By the substitution $X\mapsto {1}/{\sqrt[q]{-b}}\cdot X, Y \mapsto Y$ the polynomial $F$ goes over into a {\it normed} 
Weierstra\ss\ polynomial of type $p,q$
$$
	Y^p-X^q+\sum_{\nu p+\mu q<pq}a_{\nu\mu}X^{\nu}Y^{\mu}
$$
whose zero-set it isomorphic to $C$ and has the same place at infinity and the same Weierstra\ss\ semigroup. 
We call it the associated normed curve of $C$. 
In this paper we understand by curves of type $p,q$ the plane curves defined by normed 
Weierstra\ss\ polynomials of type $p,q$.

These curves can be identified with the points $(\{a_{\nu\mu}\}_{\nu p+\mu q<pq})\in \mathbb A^n(K)$
associated with their equation where $n:=\nfrac{1}{2}(p+1)(q+1)-1$.
In Section 1 we describe the (locally closed)
subsets of $\mathbb A^n(K)$ consisting of the various kinds of curves of type $p,q$.
In particular we are interested in the set of nodal curves of type $p,q$. Such curves have at most $d$ nodes, and 
their Weierstra\ss\ semigroup has genus $g=d-l$ if $l$ is the number of the nodes. The singular nodal curves 
form a dense open set of 
an irreducible hypersurface $\mathcal H\subset \mathbb A^n(K)$ whose properties are the main object of study in 
Section 1. It turns out that for any $l\in\{0,\dots,d\}$ a nodal curve of type $p,q$ exists 
having exactly $l$ nodes (Theorem 1.6).

Given a numerical semigroup $H$ with $p\in H$ greater than the elements of a minimal system of generators 
of $H$ we construct in Section 2 a locally closed subset $V_{pq}(H)$ in some affine
space, such that $H$ is a Weierstra\ss\ semigroup if and only if $V_{pq}(H)\ne\emptyset$ (Corollary 2.4).
The set $V_{pq}(H)$ is explicitly described by polynomial vanishing and non-vanishing conditions, 
where ''explicit'' means that the polynomials are given by a formula or there is an algorithm to 
compute them. In principle the membership test for polynomial ideals allows then to decide whether $H$ 
is a Weierstra\ss\ semigroup or not. However for any $H$ of interest (i.e.\ where the result is not known) the
number of conditions is huge so that the criterion seems only to be of theoretical interest and
not feasible for a computer program.

By the {\it simplification of nodal curves} introduced in Section 3 the criterion allows to show without computations 
that every $H$ of the following kind is a
Weierstra\ss\ semigroup: $H$ is obtained from the semigroup $H_{pq}$ generated by $p$ and
$q\ (1<p<q,$ with $p,q$ relatively prime) by closing all gaps of $H_{pq}$ which are greater than or
equal to a given gap of $H_{pq}$ (Theorem \ref{thmFamW}). The {\it hyperordinary semigroups} defined by
Rim and Vitulli [RV] belong to this class of semigroups. These authors have shown with
a different method that hyperordinary semigroups are Weierstra\ss\ semigroups.

\headA{The space of plane curves of type $p,q$}

Let $R:=K[\{A_{\nu\mu}\}_{\nu p+\mu q<pq}]$ be the polynomial ring in the $n=\nfrac{1}{2}(p+1)(q+1)-1$ 
indeterminates $A_{\nu\mu}\ (\nu p+\mu q<pq)$ over $K$. The generic (normed) Weierstra\ss\ polynomial
$$F=Y^p-X^q+\sum_{\nu p+\mu q<pq}A_{\nu\mu}X^{\nu}Y^{\mu}=A_{00}+\dots$$
of type $p,q$ has the partial derivatives
$$F_X=-qX^{q-1}+\sum_{\nu p+\mu q<pq}\nu A_{\nu\mu}X^{\nu-1}Y^{\mu}=A_{10}+\dots$$
$$F_Y=pY^{p-1}+\sum_{\nu p+\mu q<pq}\mu A_{\nu\mu}X^{\nu}Y^{\mu-1}=A_{01}+\dots $$
where the dots represent terms containing $X$ or $Y$. We are interested in the ring
$$
	A=R[X,Y]/(F,F_X,F_Y)
$$
as an $R$-Algebra. As a $K$-algebra it is isomorphic to the polynomial ring
$$
	K[\{A_{\nu\mu}\}_{(\nu,\mu)\ne (0,0),(1,0),(0,1)},X,Y]
$$
hence the image $R'$ of $R$ in $A$ is a domain. Moreover $\{F,F_X,F_Y\}$ is a regular sequence in $R[X,Y]$.

We endow $R[X,Y]$ with the grading given by deg$(X)=p$, deg$(Y)=q$ and deg$(r)=0$ for $r\in R$ and 
let $\mathcal F$ denote the corresponding degree filtration. Let $N:=R[X,Y]/(F_X,F_Y)$.
The polynomial $F$ has degree form $Y^p-X^q$, and since the partial derivatives are homogeneous 
maps the degree form of $F_X$ is $-qX^{q-1}$ and that of $F_Y$ is $pY^{p-1}$. Since they form a 
regular sequence in $R[X,Y]$ we have $\text{gr}_{\mathcal F}N=R[X,Y]/(X^{q-1},Y^{p-1})$ 
(see [Ku2], B.12), hence $N$ is a free $R$-module with the basis
$$
	B:=\{\xi^{\nu}\eta^{\mu}\}_{\nu<q-1,\mu<p-1}
$$
where $\xi,\eta$ are the residue classes of $X,Y$ in $N$ (see [Ku2], B.6). In particular $\rank N =(p-1)(q-1)=:c$, and 
different basis elements have different degrees with respect to the residue grading of $\mathcal F$.

Since $A$ is finite over $R'$ we have $R'=R/\mathfrak q$ where the prime ideal $\mathfrak q$ is 
generated by an irreducible polynomial in $R$, hence $\mathcal H:=\Spec{R'}$ is an irreducible 
hypersurface in $A^n(K)=\Spec R$.

We identify the curves of type $p,q$ with the closed points $\alpha:=(\{a_{\nu\mu}\})\in \mathbb A^n(K)$ 
or with the maximal ideals $\frak m=(\{A_{\nu\mu}-a_{\nu\mu}\}_{\nu p+\mu q<pq})\ (a_{\nu\mu}\in K)$ of $R$. 
For $\alpha\in K^n$ we denote the curve with the equation $F(\alpha,X,Y)=0$ by $C_{\alpha}$.
The set $\Max A$ can be identified with the set of singularities of the curves of type $p,q$.
If a maximal ideal $\frak M$ of $A$ with preimage $\frak m=(\{A_{\nu\mu}-a_{\nu\mu}\})$ in $R$ is given, 
then $\frak M$ corresponds to a singularity of the curve $C_{\alpha}$. Moreover
$$A_{\frak M}/\frak mA_{\frak M}=(K[C_{\alpha}]/J)_{\overline{\frak M}}$$
with the Jacobian ideal $J$ of $K[C_{\alpha}]$ and the image $\overline{\frak M}$ of $\frak M$ in $K[C_{\alpha}]/J$.

\prop{
The singular curves of type $p,q$ are the closed points of the irreducible hypersurface 
$\mathcal H\subset \mathbb A^n(K)$. The closed points of $\mathbb A^n(K)$ outside of $\mathcal H$ 
are in one-to-one correspondence with the smooth curves of type $p,q$. Their Weierstra\ss\ semigroup is $H_{pq}$.
}

The last statement of the proposition follows from the fact that the Weierstra\ss\ semigroup of a smooth curve of 
type $p,q$ has genus $g=d$, hence no gaps of $H_{pq}$ have to be closed in it.

An example of a smooth curve of type $p,q$ is given by the equation $Y^p-X^q+a_{00}=0\ (a_{00}\ne 0)$.

As an $R$-module $A$ can be written
$$
	A=N/\sum_{\alpha<q-1,\beta<p-1}R\cdot \xi^{\alpha}\eta^{\beta}F(\xi,\eta)
$$
and the relations
$$
	\xi^{q-1}=\frac{1}{q}\sum\nu A_{\nu\mu}\xi^{\nu-1}\eta^{\mu},\  
	\eta^{p-1}=-\frac{1}{p}\sum\mu A_{\nu\mu}\xi^{\nu}\eta^{\mu-1}\leqno(1)
$$
allow with the usual reduction process to write
$$
	\xi^{\alpha}\eta^{\beta}F(\xi,\eta)=\sum_{\nu<q-1,\mu<p-1}r^{\alpha\beta}_{\nu\mu}\cdot \xi^{\nu}\eta^{\mu}\qquad 
	(r^{\alpha\beta}_{\nu\mu}\in R).
$$

The  $c\times c$-matrix $M:=\left(r^{\alpha\beta}_{\nu\mu}\right)$ represents the multiplication by 
$F(\xi,\eta)$ in $N$, and since $F(\xi,\eta)$ is not a zero-divisor in $N$ we have an exact sequence of $R$-modules
$$
	0\to R^c\buildrel M \over \longrightarrow R^c\to A\to 0.\leqno (2)
$$
$M$ is a relation matrix of the $R$-module $A$ with respect to the basis $B$ of $N$. For 
$0\le l\le c$ the $(c-l)$-minors of $M$ generate the $l$-{th} Fitting ideal $F_l(A/R)$ of the $R$-module $A$. 
In particular $F_0(A/R)=(\Delta)$ with $\Delta:=\text{det}(M)$, the norm of the multiplication map by $F(\xi,\eta)$.
Here $\Delta\ne 0$, the map given by $M$ being injective. We have
$$
	(0)\ne F_0(A/R)\subset F_1(A/R)\subset \dots \subset F_c(A/R)
$$
By [Ku1], D.14
$$
	\text{Ann}_R(A)^c\subset F_0(A/R)=(\Delta)\subset \text{Ann}_R(A)\leqno (3)
$$
and therefore Rad(Ann$_R(A))=\text{Rad}(\Delta)$. As $A$ is an $R$-algebra Ann$_R(A)$ is the kernel 
$\mathfrak q$ of the structure homomorphism $R\to A$, hence Rad($\Delta)$ is also a prime ideal.
It follows that $\Delta=a\Delta_0^r$ with an irreducible polynomial $\Delta_0$ of $R$ which generates 
$\mathfrak q$, an $a\in K\setminus \{0\}$ and an $r\in \mathbb N$, hence $R'=R/\frak q=R/(\Delta_0)$ 
and the hypersurface $\mathcal H$ is given by the equation $\Delta_0=0$.

For $\frak p\in \Spec R$ and $l\in \{0,\dots,c\}$ we have the following formula for the minimal 
number of generators of the $R_{\frak p}$-module $A_{\frak p}$
$$
	\mu_{\frak p}(A)=\text{min}\{l\vert F_l(A_{\frak p}/R_{\frak p})=R_{\frak p}\}\leqno(4)
$$
([Ku1], D.8).

Let $\frak m=(\{A_{\nu\mu}-a_{\nu\mu}\})$ be a maximal ideal of $R$ corresponding to the polynomial 
$\bar F:=F(\alpha,X,Y)\in K[X,Y]$ and $l\in \{0,\dots,c\}$. Then by (4) 
$F_l(A/R)_{\frak m}=F_l(A_{\frak m}/R_{\frak m})=R_{\frak m}$ if and only if the $R_{\frak m}$-module 
$A_{\frak m}$ has a minimal number of generators $\le l$, that is, if and only if
$$
	\text{dim}_K(K[X,Y]/(\bar F,\bar F_X,\bar F_Y))\le l.\leqno(5)
$$
If $\frak M\in \Max{K[X,Y]}$ corresponds to a node of $C_{\alpha}$, then
$$\text{dim}_K((K[X,Y]/(\bar F,\bar F_X,\bar F_Y))_{\frak M})=1.$$
If $C_{\alpha}$ is a nodal curve, then dim$_K(K[X,Y]/(\bar F,\bar F_X,\bar F_Y)$) is the number of its nodes
and (4) implies

\lma{\label{LmaAtMost}
If $C_{\alpha}$ has at most $l$ nodes and no other singularities, then $\frak m$ 
is contained in the open set $\Max R \setminus V(F_l(A/R))$ of  $\Max R$. Conversely, if $C_{\alpha}$ has 
$l$ distinct nodes and $\mathfrak m\in\Max R\setminus V(F_l(A/R))$, then $C_{\alpha}$ is a nodal 
curve with exactly $l$ nodes.
}

For the module of differentials of $A/R$ we have
$$
	\Omega^1_{A/R} = AdX\oplus AdY/ \left\langle F_{XX}(x,y)dX+F_{XY}(x,y)dY, F_{YX}(x,y)dX+F_{YY}(x,y)dY\right\rangle
$$
with the residue classes $x,y$ of $X,Y$ in $A$.
{Since the variables $A_{00}, A_{10}, A_{01}$ have disappeared in the second derivatives the Hesse determinant 
$\Hess F(X,Y)$ of $F$ is a non-zero polynomial in $A$.}
Take $\frak M\in \Max A$ with preimage $\frak m$ in $R$ corresponding to a point in $\mathcal H$. 
Nodes are the singularties with non-vanishing Hesse determinant, hence $\frak M$ corresponds to a node of 
the curve given by $\frak m$ if and only if $\frak M\in \Max A\setminus V(\Hess F)$. 
This is equivalent to each of the following conditions
\begin{enumerate}
	\item[(i)] $\Hess F(X,Y)$ is a unit in $A_{\frak M}$.
	\item[(ii)] $\Omega^1_{A_{\mathfrak M}/R}=0$.
	\item[(iii)] $\frak M$ is unramified over $R$ ([Ku1], 6.10).
\end{enumerate}

From (ii) and (iii) we conclude

\prop{
The nodal curves of type $p,q$ having at least one node correspond bijectively to the maximal 
ideals $\mathfrak m\in V(\Delta_0)=\mathcal H$ with $\frak m\not\in V(\text{Ann}_R(\Omega^1_{A/R}))$ 
or equivalently with $A/R$ being unramified at $\mathfrak m$.
}

We denote this open set of the hypersurface $\mathcal H$ by $\mathcal H_{pq}$.
The additional assumption that $\frak m\not\in V(F_l(A/R))$ defines for each $l=1,\dots,d$
an open subset $U_l$ of $\mathcal H_{pq}$ whose closed points correspond to the nodal curves of type $p,q$ 
having at most $l$ nodes. Set $U_0:=\emptyset$. We have 
$$
	\mathcal H_{pq}=\bigcup_{l=1}^d\mathcal H^l_{pq}
$$
with the locally closed subset $\mathcal H^l_{pq}:= U_l\setminus U_{l-1}$ whose closed points correspond 
to the curves having exactly $l$ nodes. The Weierstra\ss\ semigroups of the curves in $\mathcal H^l_{pq}$ are certain 
semigroups $H$ with $p,q\in H$ having genus $g=d-l$. By [KKW], Theorem 6.4 every Weierstra\ss\ semigroup $H$ 
of genus $g$ is the Weierstra\ss\ semigroup of an element of $\mathcal H_{pq}^{d-g}$ for suitably chosen $p,q$.

The curve $C: (Y-b)^p-(X-a)^q+c(X-a)(Y-b)=0\ (a,b,c\in K,c\ne 0))$ has only one singularity at $(a,b)$, 
and it is a node. Therefore $\mathcal H^1_{pq}\ne \emptyset$. The associated normed curve of the 
Lissajous curve of type $p,q$ ([KKW], Example 2.4) has the maximal possible number $d$ of nodes, hence 
$\mathcal H^d_{pq}\ne \emptyset$.

For a domain $B$ let $Q(B)$ denote its quotient field.

\prop{
We have $Q(R')=Q(A)$. Hence $R'\to A$ induces a finite birational morphism 
$\mathbb A^{n-1}(K)\to \mathcal H$, and the hypersurface $\mathcal H$ is rational.
}

\prf{
Let $\frak p$ be the kernel of $R\to A$. The inclusion $R'\to A$ induces an injection $Q(R')\to A_{\frak p}$. 
Since $A$ is a domain and integral over $R'$ we have $A_{\frak p}=Q(A)$. Moreover 
$F_1(A_{\frak p}/R_{\frak p})=R_{\frak p}$ since $F_1(A_{\frak m}/R_{\frak m})=R_{\frak m}$ with the 
$\frak m$ belonging to the curve $C$ above, as Fitting ideals are compatible with localization. 
Hence by (4) $A_{\frak p}$ is generated over $R_{\frak p}$ by one element, i.e.\ $Q(A)=Q(R')$.
}\mbox{}\\

For the maximal ideals $\frak m\in \mathcal H^l_{pq}$ all $\frak M\in \Max A$ lying over $\frak m$ are unramified 
over $R$. 
Let $\mathfrak m=(\{A_{\nu\mu}-a_{\nu\mu}\}_{\nu p+\mu q<pq})\in \mathcal H^l_{pq}$ with 
$\alpha:=(\{a_{\nu\mu}\})\in K^n$ be given, and let $\frak M\in \Max A$ correspond to a node of 
the curve $C_{\alpha}$. Set $T_{\nu\mu}:=A_{\nu\mu}-a_{\nu\mu}$ for short.
The canonical homomorphism $R_{\mathfrak m}\to A_{\mathfrak M}$ induces a local homomorphism 
$\varphi: \widehat{R_{\mathfrak m}}\to \widehat{A_{\mathfrak M}}$ of the completions which is 
surjective since $A_{\mathfrak M}$ is unramified over $R_{\mathfrak m}$. Here 
$\widehat{R_{\mathfrak m}}=K[[\{T_{\nu\mu}\}]]$ and $\widehat{A_{\mathfrak M}}$ are regular 
local rings of dimension $n$ resp.\ $n-1$. Therefore ker$(\varphi)$ is generated by a power 
series $\Delta_{\mathfrak M}$ of order 1, an irreducible factor of $\Delta_0$ considered as a power series 
in the $T_{\nu\mu}$. Thus $\mathfrak M$ defines a smooth analytic branch 
$\Spec{\widehat{R_{\mathfrak m}}/(\Delta_{\mathfrak M})}$ of $\mathcal H$ near the point $\alpha$. 
Different nodes of $C_{\alpha}$ define different branches as the power series $\Delta_0$ cannot 
have multiple factors, $\Delta_0$ being an irreducible polynomial. The 
local ring $R_{\mathfrak m}/(\Delta_0)$ is regular if and only if $C_{\alpha}$ has only one node. 
Thus $\mathcal H^1_{pq}$ is the set of regular points of $\mathcal H_{pq}$.
\vskip3mm
Let $\widehat{A_{\mathfrak m}}$ be the completion of $A_{\mathfrak m}:=R_{\mathfrak m}\otimes A$ 
as an $R_{\mathfrak m}$-module. Then 
$$
	\widehat{A_{\mathfrak m}}=\widehat{A_{\mathfrak M_1}}\times\cdots\times\widehat{A_{\mathfrak M_l}}=
		\widehat{R_{\mathfrak m}}/(\Delta_{\mathfrak M_1})\times\cdots\times
		\widehat{R_{\mathfrak m}}/(\Delta_{\mathfrak M_l})\leqno(6)
$$
by the Chinese Remainder Theorem.
Since the Fitting ideals are compatible with localization and completion we obtain that 
$$
	F_0(\widehat{A_{\mathfrak m}}/\widehat{R_{\mathfrak m}})=\widehat{R_{\mathfrak m}}\cdot F_0(A/R)=
			\widehat{R_{\mathfrak m}}\cdot\Delta=\widehat{R_{\mathfrak m}}\cdot(\prod_{i=1}^l\Delta_{\mathfrak M_i})=
			\widehat{R_{\mathfrak m}}\cdot\Delta_0.
$$
Remember that $\Delta=a\Delta_0^r$ with $a\in K\setminus\{0\}$ and $r\ge 1$. Since we know that nodal curves 
of type $p,q$ with at least one node exist for every $p,q$ the above consideration implies that $r=1$ and we have 
proved the irreducibility of $\Delta$, the polynomial generating $F_0(A/R)$. Thus the hypersurface 
$\mathcal H$ is defined by $F_0(A/R)=(\Delta)$.

We determine the leading form of $\Delta_{\mathfrak M}$ which defines the tangent hyperplane of the branch 
$\Delta_{\mathfrak M}=0$. As $\Omega^1_{A_{\mathfrak M}/R}=0$ we see that $\Omega^1_{A_{\mathfrak M}/K}$ is 
generated by the differentials $dT_{\nu\mu}\ (\nu p+\mu q<pq)$. Moreover we have the relation
$$
	\sum_{\nu p+\mu q<pq}x^{\nu}y^{\mu}dT_{\nu\mu}=0 \leqno(7)
$$
coming from $dF=0$. Therefore $\{dT_{\nu\mu}\}_{(\nu,\mu)\ne (0,0)}$ is a basis of $\Omega^1_{A_{\mathfrak M}/K}$, 
and we obtain
$$
	\widehat{\Omega^1_{A_{\mathfrak M}/K}}=\bigoplus_{(\nu,\mu)\ne (0,0)}\widehat{A_{\mathfrak M}}\cdot dT_{\nu\mu}
$$
In $\widehat{\Omega^1_{A_{\mathfrak M}/K}}$ there is the relation 
$d\Delta_{\mathfrak M}=\sum_{\nu p+\mu q<pq}\partial \Delta_{\mathfrak M}/\partial T_{\nu\mu}\cdot dT_{\nu\mu}$,
and by (7)
$$
	\sum_{(\nu,\mu)\ne (0,0)}\bigl(\partial \Delta_{\mathfrak M}/\partial T_{\nu\mu}-x^{\nu}y^{\mu}\cdot
	\partial\Delta_{\mathfrak M}/\partial T_{00}\bigr)dT_{\nu\mu}=0
$$
which implies in $\widehat{A_{\mathfrak M}}$ the relations
$$
	\partial\Delta_{\mathfrak M}/\partial T_{\nu\mu}=x^{\nu}y^{\mu}\cdot\partial\Delta_{\mathfrak M}/\partial T_{00}
$$
for all $\nu,\mu$. Since $\Delta_{\mathfrak M}$ has order 1, at least one of the partial derivatives must be a 
unit in $\widehat{A_{\mathfrak M}}$, hence so must be the partial with respect to $T_{00}$. 
Let $(\xi,\eta)\in K^2$ be the node corresponding to $\mathfrak M$. Considering the above relations modulo 
$\mathfrak M\widehat{A_{\mathfrak M}}$ we find that
$$
	\partial\Delta_{\mathfrak M}/\partial T_{\nu\mu}\vert_0=\xi^{\nu}\eta^{\mu}
	\cdot\partial\Delta_{\mathfrak M}/\partial T_{00}\vert_0
$$ where the last partial does not vanish,
hence the leading form of $\Delta_{\mathfrak M}$ is 
$$
	L_{\mathfrak M}\Delta_{\mathfrak M}=\partial\Delta_{\mathfrak M}/\partial T_{00}\vert_0\cdot 
	\sum_{\nu p+\mu q<pq}\xi^{\nu}\eta^{\mu}T_{\nu\mu}.\leqno(8)
$$
Collecting everything we obtain

\prop{
At the closed points of $\mathcal H^l_{pq}$ the hypersurface $\mathcal H$ has $l$ regular branches 
with tangent hyperplanes given by (8).
}

Let $(\xi_1,\eta_1),\dots,(\xi_l,\eta_l)$ be the nodes of $C_{\alpha}$. 
We shall see in Lemma 2.1 that the matrix 
$\left(\xi_i^{\nu}\eta_i^{\mu}\right)_{\nu p+\mu q<pq, i=1,\dots,l}$ has rank $l$. 
Thus the $L_{\mathfrak M_i}\Delta_{\mathfrak M_i}\ (i=1,\dots,l)$ are linearly independent over 
$K$ and the $\Delta_{\mathfrak M_i}$ form part of a regular system of parameters of $\widehat{R_{\mathfrak m}}$.

For the defining polynomial $\Delta$ of $\mathcal H$ this means the following: If we expand $\Delta$ as a 
polynomial in the $T_{\nu\mu}=A_{\nu\mu}-a_{\nu\mu}$ its form of lowest degree is  up to a constant factor 
the product of the $l$ homogenous linear polynomials $L_{\mathfrak M_i}\Delta_{\mathfrak M_i}$ which 
moreover are linearly independent over $K$.

Formula (6) implies that
$$
	\widehat{R_{\mathfrak m}}\cdot F_k(A/R)=F_k(\widehat{A_{\mathfrak m}}/\widehat{R_{\mathfrak m}})=
	(\{\Delta_{\mathfrak M_{i_1}}\cdot \dots \cdot \Delta_{\mathfrak M_{i_{l-k}}}\}_{i_1<\dots <i_{l-k}}).
$$

Let $\mathfrak p_{i_1,\dots,i_k}$ be the prime ideal of $\widehat{R_{\mathfrak m}}$ generated by 
$\Delta_{\mathfrak M_{i_1}},\dots,\Delta_{\mathfrak M_{i_k}}$ where $i_1<\dots <i_k$. Then
$$
	\widehat{R_{\mathfrak m}}\cdot F_k(A_{\mathfrak m}/R_{\mathfrak m})= 
	\bigcap_{i_1<\dots <i_{k+1}}\mathfrak p_{i_1,\dots,i_{k+1}}\ (k=0,\dots,l-1)\leqno(9)
$$
in particular
$$
	\widehat{R_{\mathfrak m}}\cdot F_{l-1}(A/R)=F_{l-1}(\widehat{A_{\mathfrak m}}/\widehat{R_{\mathfrak m}})
		=(\Delta_{\mathfrak M_1},\dots,\Delta_{\mathfrak M_l})=\mathfrak p_{1,\dots,l}.
$$
One can prove (9) by first showing it when the $\Delta_{\mathfrak M_i}$ are variables in a polynomial ring and  
by passing then to the completion.
Thus the ideals $\widehat{R_{\mathfrak m}}\cdot F_k(A/R)\ (k=0,\dots,l-1)$ are radical ideals of height 
$k+1$ in $\widehat{R_{\mathfrak m}}$, and so are the  $F_k(A_{\mathfrak m}/R_{\mathfrak m})$ in $R_{\mathfrak m}$.

\thm{
For any $l$ with $1\le l\le d$ there is a nodal curve of type $p,q$ with exactly $l$ nodes, i.e.\ 
$\mathcal H^l_{pq}\ne \emptyset$.
}

\prf{
Let $\mathfrak m\in \Max R$ correspond to the curve $C$ associated to the Lissajous curve of type $p,q$. 
There is a $g\in R$ such that $\mathfrak m\in D(g)$ and that the closed points in $D(g)\is \mathcal H$ correspond to nodal curves. 
Then by the above the $F_k(A_g/R_g)\ (k=0,\dots,d)$ form a strictly increasing sequence of radical ideals in $R_g$. 
Choose a maximal ideal $\mathfrak n\in D(g)$ such that
$$F_{l-1}(A_g/R_g)\subset \mathfrak n R_g,\ F_l(A_g/R_g)\not\subset \mathfrak n R_g.$$
Then the curve corresponding to $\mathfrak n$ has exactly $l$ nodes. 
}

\expls{
The Weierstra\ss\ semigroups of the curves in $\mathcal H^l_{pq}$ are numerical semigroups 
$H$ with $p,q\in H$ and genus $d-l$. If $l\le \nfrac{p}{2}$ all possible $H$ of this kind do occur, see [KKW], Example 5.4. 
For $l=1$ the semigroup $H$ is obtained from
$H_{pq}$ by closing one gap $c-1-(ap+bq), (a,b\in \mathbb N)$. We must have $a=b=0$, otherwise more than one gap would be closed. 
Therefore $H=\left\langle p,q,c-1 \right\rangle$ and any curve in $\mathcal H^1_{pq}$ has this Weierstra\ss\ semigroup.  In $\mathcal H^2_{p,q}$ 
we have the Weierstra\ss\ semigroups $\langle p,q,c-1-p\rangle$ and $\langle p,q,c-1-q\rangle$.
The curves in $\mathcal H^d_{pq}$ are the nodal curves of type $p,q$ for which the normalization of its projective 
closure has genus $0$. Their Weierstra\ss\ semigroup is $\mathbb N$.
}

The hypersurface $\mathcal H$ contains many lines.
\prop{
Let $\alpha \not= \beta$ be closed points of $\mathcal H$ such that 
$\Sing{C_\alpha} \is \Sing{C_\beta} \not= \emptyset$, 
and let $L$ be the line through $\alpha$ and $\beta$. Then $L \subset \mathcal H$, and for
almost all closed $\gamma \in L$ the curve $C_\gamma$ has the singular set $\Sing{C_\alpha} \is \Sing{C_\beta}$.
}

\prf{
Set $H := F(\beta,X,Y) - F(\alpha,X,Y)$ and $D := V (H)$. Then $H$ and
$F(\alpha,X,Y)$ are relatively prime and $\Sing{C_\alpha} \is \Sing{D} = \Sing{C_\alpha} \is \Sing{C_\beta}$.
By [KKW], Proposition 3.1 the curve
\[
	F(\alpha,X,Y) + d \cdot H = F(\alpha + d(\beta - \alpha),X,Y) = 0
\]
has for almost all $d \in K \setminus\{0\}$ the singular set $\Sing{C_\alpha} \is \Sing{C_\beta} \not= \emptyset$. It
follows that $L \subset \mathcal H$.
}

\cor{
For any closed point $\alpha \in \mathcal H$ there is at least one line $L$ with $\alpha \in L \subset \mathcal H$.
}

\prf{
Let $(a,b)$ be a singularity of $C_\alpha$, and let $C_\beta$ be a nodal curve with $(a,b)$
as its only node. It can be chosen such that $\alpha \not= \beta$. Then $H$ contains by 1.8 the
line through $\alpha$ and $\beta$.
}

\cor{
Let $L \subset \mathcal H$ be a line through a closed point $\alpha$ where $C_\alpha$ is a
nodal curve. Then for almost all closed points $\gamma \in L$ the curves $C_\gamma$ have the
same Weierstra{\ss} semigroup.
}

\prf{
Since $\alpha \in \mathcal H_{pq}$ and this set is open in $\mathcal H$ almost all $C_\gamma$ with $\gamma \in L$ are
nodal curves having by 1.8 the same set of nodes. By [KKW], Corollary 4.3
they also have the same Weierstra{\ss} semigroup.
}

\headA{Which numerical semigroups are Weierstra\ss semi- \\groups?}

Let $H$ be a numerical semigroup of genus $g$ and let $p<q$ be relatively prime numbers from $H$.
The semigroup $H_{pq}$ has $d$ gaps $\gamma_1<\dots <\gamma_d$ which can be written
$$\gamma_i=(p-1)(q-1)-1-(a_ip+b_iq)\ (i=1,\dots,d)$$
with a unique $(a_i,b_i)\in \mathbb N^2.$ Of these gaps $l:=d-g$ are closed in $H$. We want to decide 
whether a nodal curve $C$ of type $p,q$ with $l$ nodes exists such that $H$ is the Weierstra\ss\ semigroup of $C$.

Let $\gamma_{j_1} < \cdots < \gamma_{j_l}$ be the gaps of $H_{pq}$ which are closed in $H$.
Further let \linebreak $\mathcal A_H(X_1,Y_1,\dots,X_l,Y_l)$ be the matrix
$$
	\left( \begin{array}{ccc}
			X_1^{a_{j_1}}Y_1^{b_{j_1}}&\dots& X_1^{a_{j_l}}Y_1^{b_{j_l}}\\
			\vdots	& \cdots& \vdots\\
			X_l^{a_{j_1}}Y_l^{b_{j_1}}&\dots&X_l^{a_{j_l}}Y_l^{b_{j_l}}
	\end{array}\right)
$$
and $D_H(X_1,Y_1,\dots,X_l,Y_l):=\det{\mathcal A_H(X_1,Y_1,\dots,X_l,Y_l)}$
its determinant.

\lma{
If $H$ is the Weierstra\ss\ semigroup of a nodal curve $C: F=0$ of type $p,q$ 
with the nodes $(\xi_1,\eta_1),\dots,(\xi_l,\eta_l)$, then $D_H(\xi_1,\eta_1,\dots,\xi_l,\eta_l)\ne 0$.
}

\prf{
If $D_H(\xi_1,\eta_1,\dots,\xi_l,\eta_l)=0$, then there exists a non-zero 
$\lambda=(\lambda_1,\dots,\lambda_l)$ $\in K^l$ 
such that $\mathcal A_H\cdot \lambda^t=0$. Assume that $\lambda_1=\cdots =\lambda_{i-1}=0,\ \lambda_i\ne 0$. 
Let $x, y$ denote the images of $X, Y$ in the function field $K(C)$ of $C$ and $P$ the place at
infinity of $C$.
The function
$$
	\Phi(x,y):=\lambda_ix^{a_{j_i}}y^{b_{j_i}}+\dots
		+\lambda_lx^{a_{j_l}}y^{b_{j_l}}\in K[C]
$$
satisfies $\Phi(\xi_i,\eta_i)=0\ (i=1,\dots,l)$. If follows from [KKW], Proposition 4.2 that \linebreak 
$\ordP{\frac{\Phi(x,y)}{F_Y(x,y)}dx} + 1 = \gamma_{j_i}$ is a gap of $H$,
contradicting the fact that $\gamma_{j_i}$ was a gap of $H_{pq}$ closed in $H$.
}\mbox{}\\

With the generic Weierstra\ss\ polynomial $F(\{A_{\nu\mu}\},X,Y)\in R[X,Y]$ of type $p,q$ and $l$ with 
$1\le l\le d$ set 
$$
	T:=R[X_1,Y_1,\dots,X_l,Y_l]/(\{F(X_i,Y_i),F_X(X_i,Y_i),F_Y(X_i,Y_i)\}_{i=1,\dots,l}).
$$

Let $C_L: F(\{a^L_{\nu\mu}\},X,Y)=0$ be the normed curve associated to the Lissajous curve of 
type $p,q$, and let  $(\xi_i,\eta_i)\ (i=1,\dots,d)$ be its nodes. $C_L$ has the 
Weierstra\ss\ semigroup $\mathbb N$. By Lemma 2.1 we have $D_{\mathbb N}(\xi_1,\eta_1\dots,\xi_d,\eta_d)\ne 0$. 
Therefore the columns of this determinant corresponding to the gaps $\gamma_{j_1},\dots,\gamma_{j_l}$ 
are linearly independent over $K$, and there are nodes $(\xi_1,\eta_1)\dots,(\xi_l,\eta_l)$ (say) such that 
$D_H(\xi_1,\eta_1,\dots,\xi_l,\eta_l)\ne 0$ too.

Let $\delta$ be the image of $D_H(X_1,Y_1,\ldots,X_l,Y_l)$ and $t$ that of 
$\prod_{i=1}^l\text{Hess}_F(X_i,Y_i)$ in $T$. Then $t\cdot \delta$ is not contained 
in the maximal ideal corresponding to the point $(\{a_{\nu\mu}^L\},\xi_1,\eta_1,\dots,\xi_l,\eta_l)$ 
and hence $t\cdot \delta$ is not nilpotent.
Therefore $S_{\!H}:=T_{t\cdot\delta}$ is not the zero-ring. Now the elements of $\Max{S_{\!H}}$ 
correspond bijectively to the $(\beta,\xi_1,\eta_1,\dots,\xi_l,\eta_l)$ where the $(\xi_i,\eta_i)$ 
are nodes of the curve $C_{\beta}$ and have the additional property that 
$D_H(\xi_1,\eta_1,\dots,\xi_l,\eta_l)\ne 0$. In particular the nodes are distinct.\\

Let $h$ be the set of the $(a_{j_i},b_{j_i})\in \mathbb N^2\ (i=1,\dots,l)$ corresponding to the gaps of $H_{pq}$ 
which are closed in $H$. 
Let $x_i,y_i$ be the images of the $X_i,Y_i$ in $S_{\!H}$ and denote the images of the $A_{\nu\mu}$ also by 
$A_{\nu\mu}\ (\nu p+\mu q<pq)$.

\lma{
We have $\Omega^1_{S\!_{H}/K}=\bigoplus_{(\nu,\mu)\not\in h}S_{\!H} dA_{\nu\mu}$. In particular $S_{\!H}$ is a regular 
$K$-algebra, equidimensional of dimension $n-l$. Further $S_{\!H}$ is unramified over 
$K[\{A_{\nu\mu}\}_{(\nu,\mu)\not \in h}]$.
}

\prf{
The module of differentials has the presentation
$$
	\Omega^1_{S\!_{H}/K}=\bigoplus_{\nu p+\mu q<pq}S_{\!H}dA_{\nu\mu}\oplus\bigoplus_{i=1}^lS_{\!H}dX_i\oplus S_{\!H}dY_i/U
$$
where U is generated by
$$
	\sum_{\nu p+\nu q<pq}x_i^{\nu}y_i^{\mu}dA_{\nu\mu},
$$
$$
	\sum_{\nu p+\nu q<pq}\nu x_i^{\nu-1}y_i^{\mu}dA_{\nu\mu}+F_{XX}(x_i,y_i)dX_i+F_{XY}(x_i,y_i)dY_i
$$
and
$$
	\sum_{\nu p+\nu q<pq}\mu x_i^{\nu}y_i^{\mu-1}dA_{\nu\mu}+F_{YX}(x_i,y_i)dX_i+F_{YY}(x_i,y_i)dY_i
$$
($i=1,\dots,l$). Since $\Hess F(x_i,y_i)\ (i=1,\dots,l)$ and $D_H(x_1,y_1,\dots,x_l,y_l)$ 
are units in $S_{\!H}$ the statement about $\Omega^1_{S\!_{H}/K}$ follows, and the remaining assertions are clear by the 
differential criterion of regularity ([Ku1],7.2).
}\mbox{}\\

Now let $U^l_{pq}(H):=\Spec {S_{\!H}}\setminus V(F_l(A/R)S_{\!H})$. By Lemma \ref{LmaAtMost} the closed points \linebreak
$(\alpha,\xi_1,\eta_1,\dots,\xi_l,\eta_l)$ of the scheme $U^l_{pq}(H)$ are those for which the curve 
$C_{\alpha}$ has no singularities but the nodes $(\xi_i,\eta_i)$ which satisfy 
$D_H(\xi_1,\eta_1,\dots,\xi_l,\eta_l)\ne 0$. 
These $C_{\alpha}$ have a Weierstra\ss\ semigroup which is obtained from $H_{pq}$ by closing $l$ of its gaps, 
but may be different from $H$.

It will be shown in Proposition \ref{propU} that the scheme $U^l_{pq}(H)$ is not empty. 
In order to decide whether $H$ is the Weierstra\ss\ semigroup of a nodal curve of type $p,q$ we need a further 
consideration which is inspired by [Ha], IV.4.

Let $\gamma_{i_1}<\dots<\gamma_{i_g}$ be the gaps of $H,\ \gamma_{i_k}=c-1-(a_{i_k}p+b_{i_k}q)$. Then
$$
	\{\gamma_{i_1},\dots,\gamma_{i_g}\}\cup\{\gamma_{j_1}\dots,\gamma_{j_l}\}
$$
is the set of all gaps of $H_{pq}$.
In $H_{pq}$ there are $d-i_k$ gaps $>\gamma_{i_k}$, and $H$ has $g-k$ gaps $>\gamma_{i_k}$. Hence there are 
$(d-i_k)-(g-k)=l-(i_k-k)$ gaps of $H_{pq}$ which are $>\gamma_{i_k}$ and are closed in $H$. Therefore 
$\gamma_{j_m}>\gamma_{i_k}$ if and only if $m>i_k-k$.

Let $s_k$ be the column
$$
	\left( \begin{array}{c}
		X_1^{a_{i_k}}Y_1^{b_{i_k}}\\
		\vdots\\
		\vdots\\
		X_l^{a_{i_k}}Y_l^{b_{i_k}}
	\end{array}\right) \ (k=1,\dots,g)
$$
and $D^m_k(X_1,Y_1,\dots,X_l,Y_l)$ for $m\in \{1,\dots,i_k-k\}$ the determinant of the matrix which is obtained from 
$\mathcal A_H$ by replacing its m-th column by $s_k$. These are 
$\sum_{k=1}^g (i_k-k)=\sum_{k=1}^g i_k-{g+1\choose 2}$ determinants. 
Let $J$ be the ideal generated by their images in $S_{\!H}$. If the semigroup $H$ is obtained from 
$H_{pq}$ by closing its $l$ greatest gaps, then no $D^m_k$ are present, and we set $J=(0)$.
Let $V_{pq}(H):=U^l_{pq}(H)\cap V(J)$.

\thm{
The closed points of $V_{pq}(H)$ correspond to the nodal curves of type $p,q$ 
having the Weierstra\ss\ semigroup $H$, i.e.\ $H$ is the Weierstra\ss\ semigroup of such a curve if and 
only if $V_{pq}(H)\ne \emptyset$.
}

\prf{
a) Let $Q:=(\alpha,\xi_1,\eta_1,\dots,\xi_l,\eta_l)\in V_{pq}(H)$. Since $Q\in U^l_{pq}(H)$ the curve $C_{\alpha}$ 
is a nodal curve with the nodes $(\xi_1,\eta_1),\dots,(\xi_l,\eta_l)$ and 
$D_H(\xi_1,\eta_1,\dots,\xi_l,\eta_l)\ne 0$. Moreover
$$
	D^m_k(\xi_1,\eta_1,\dots,\xi_l,\eta_l)=0\ \text{for}\ k=1,\dots,g\ \text{and}\ m=1,\dots,i_k-k.\leqno(1)
$$

Further for any $k\in\{1,\dots,g\}$ the linear system of equations
$$
\mathcal A_H(\xi_1,\eta_1,\dots,\xi_l,\eta_l)
\left( \begin{array}{c} 
	\lambda_1\\
	\vdots\\
	\lambda_l
\end{array}\right) = -s_k(\xi_1,\eta_1,\dots,\xi_l,\eta_l)
$$
has a unique solution. By Cramer's rule (1) implies
that $\lambda_1=\cdots =\lambda_{i_k-k}=0$. Let $x,y$ denote the images of $X,Y$ in the function field of $C_{\alpha}$.
The polynomial
$$
	\Phi_k(X,Y):= X^{a_{i_k}}Y^{b_{i_k}}+\sum_{m>i_k-k}\lambda_mX^{a_{j_m}}Y^{b_{j_m}}
$$
vanishes at the nodes $(\xi_i,\eta_i)\ (i=1,\dots,l)$,
and since $\gamma_{j_m}>\gamma_{i_k}$ for $m> i_k-k$ the differential $\omega_k:=\frac{\Phi_k(x,y)}{F_Y(x,y)}dx$ 
has order $\ordP{\omega_k} = \gamma_{i_k}-1$ at the place at infinity of $C_{\alpha}$. 
By [KKW], Proposition 4.2 $\gamma_{i_1},\dots,\gamma_{i_g}$ are gaps of the Weierstra\ss\ semigroup of 
$C_{\alpha}$, i.e.\ $H$ is this semigroup.\\

b) Let $H$ be the Weierstra\ss\ semigroup of a nodal curve $C_{\alpha}: F(\alpha,X,Y)=0$ of type $p,q$ with 
$l$ distinct nodes $(\xi_1,\eta_1),\dots,(\xi_l,\eta_l)$. We show that 
$Q:=(\alpha,\xi_1,\eta_1,\dots,\xi_l,\eta_l)\in V_{pq}(H)$. 
By the discussion above we know already that $Q\in U^l_{pq}(H)$.

Let $\Omega_{\infty}$ be the vector space of differentials with non-negative order at the place $P$ at infinity of 
$C_{\alpha}$. According to [KKW], Lemma 4.1 we can choose
a basis $\{\omega_1,\dots,\omega_l\}$ of the vector space $\Omega$ of holomorphic differentials on $\mathcal R$
such that 
$\omega_k=\frac{\Phi_k(x,y)}{F_Y(x,y)}dx$ with
$$	
	\Phi_k(x,y)=x^{a_{i_k}}y^{b_{i_k}} + \lambda_{i_k+1}x^{a_{i_k+1}}y^{b_{i_k+1}}+\dots +
	\lambda_dx^{a_d}y^{b_d}\ (k=1,\dots,g)
$$
and $\ordP{\omega_k} + 1 = \gamma_{i_k}$.
By elementary transformations we attain that
$$
	\Phi_k(x,y)=x^{a_{i_k}}y^{b_{i_k}}+\tilde\lambda_{r,k}x^{a_{j_r}}y^{b_{j_r}}+\dots
	+\tilde\lambda_{l,k}x^{a_{j_l}}y^{b_{j_l}}\ (k=1,\dots,g),
$$
with certain $\tilde\lambda_{i,k}\in K$
where $r=i_k-k+1$. Since $\Phi_k(\xi_i,\eta_i)=0\ (i=1,\dots,l)$ and $\tilde\lambda_{k,m}=0\ (m=1,\dots,i_k-k)$
Cramer's rule implies that $D^m_k(\xi_1,\eta_1,\dots,\xi_l,\eta_l)=0$ for $k=1,\dots,g$ and $m=1,\dots,i_k-k$. 
Hence $Q\in V(J)\cap U^l_{pq}(H)=V_{pq}(H)$.
}\mbox{}\\

Theorem 2.3 and [KKW], Theorem 6.4 imply

\cor{
Let $p$ be greater than the elements of the minimal system of generators of $H$. 
Then $H$ is a Weierstra\ss\ semigroup if and only if $V_{pq}(H)\ne \emptyset$.
}

The closed points of $V_{pq}(H)$ are the 
$(\{a_{\nu\mu}\}_{\nu p+\mu q<pq},\xi_1,\eta_1,\dots,\xi_l,\eta_l)\in K^{n+2l}$
which are zeros of the polynomials
$$
	F(X_i,Y_i),\ F_X(X_i,Y_i),\ F_Y(X_i,Y_i)\	(i=1,\dots,l)\leqno(2)
$$
and of
$$
	D^m_k(X_1,Y_1,\dots,X_l,Y_l)\ (k=1,\dots,g,m=1,\dots,i_k-k)\leqno(3)
$$
and not zeros of the polynomials $D_H(X_1,Y_1,\dots,X_l,Y_l),
\Hess F(X_i,Y_i)\ (i=1,\dots,l)$
and of at least one of the $N:={c\choose l}^2\ (c-l)$-minors $h_t\ (t=1,\dots,N)$ of the matrix 
$M= \left(r^{\alpha\beta}_{\nu\mu}\right)$ defined in Section 1.
Let $\mathfrak a$ be the ideal in $R[X_1,Y_1,\dots,X_l,Y_l]$ generated by the polynomials (2) and (3).
By Theorem 2.3 and Hilbert's Nullstellensatz $H$ is the Weierstra\ss\ semigroup of a nodal curve of type $p,q$
if and only if there exists $t\in \{1,\dots,N\}$ such that
$$
	h_t\cdot\prod_{i=1}^l\Hess F(X_i,Y_i)\cdot D_H(X_1,Y_1,\dots,X_l,Y_l)\not\in \Rad {\mathfrak a}.\leqno(4)
$$
One can try to decide this by the radical membership test (see e.g. [Kr-R], page 219). 
However the number $N$ of necessary tests increases rapidly with $p$ and $q$, 
and so do the degrees of the involved polynomials. 
A sufficient condition is that (4) holds for a $(c-1)$-minor $h_t$ of the matrix $M$ which requires $c^2$ tests 
in the worst case, but with no guarantee of a success.

The polynomials in (2),(3) and (4) all belong to $\mathbb Q[\{A_{\nu\mu}\},X_1,Y_1,\dots,X_l,Y_l]$. 
Therefore (4) holds true if and only if it holds true for $K=\overline{\mathbb Q}$, the field of algebraic numbers. 
In other words, the property of $H$ to be a Weierstra\ss\ semigroup is independent of the choice of the base field. 
For example we can test it for $K=\mathbb C$.

The projection $\mathbb A^{n+2l}(K)\to \mathbb A^n(K)\ ((\alpha,\xi_1,\eta_1,\dots,\xi_l,\eta_l)\mapsto \alpha)$ 
maps the locally closed set $V_{pq}(H)$ onto a {\it constructible} set $V^H_{pq}\subset \mathcal H^l_{pq}$ whose 
closed points correspond bijectively to the nodal curves of type $p,q$ with the Weierstra\ss\ semigroup $H$. We have
$$
	\mathcal H^l_{pq}=\bigcup_HV^H_{pq}
$$
where $H$ runs over the numerical semigroups containing $p$ and $q$ with $d-l$ gaps.

\headA{Simplification of nodal curves and a class of Weierstra\ss\ semigroups}

Let $1<p<q$ be relatively prime integers and $d= \nfrac{1}{2}(p-1)(q-1)$. In Theorem 1.6 we have seen that for any 
$l\in \{1,\dots,d\}$ there is a nodal curve of type $p,q$ with exactly $l$ nodes. 
The following proposition gives a more precise statement and a different proof.

\prop{\label{propU}
Let $H$ be a numerical semigroup which is obtained from $H_{pq}$ by closing $l$ of its gaps. 
Then $U^l_{pq}(H)\ne \emptyset$.
}

As an immediate consequence we get

\thm{\label{thmFamW}
Let $H$ be the numerical semigroup which is obtained from $H_{pq}$ by closing its $l$ 
greatest gaps. Then $H$ is a Weierstra\ss\ semigroup.
}

{
In fact, for $H$ as in \ref{thmFamW} no determinants $D^m_k$ occur. Therefore $V_{pq}(H)=U^l_{pq}(H)$ which is not empty 
by \ref{propU}, and Theorem 2.3 implies that $H$ is the Weierstra\ss\ semigroup of a nodal curve of type $p,q$.
}\hfill \ensuremath{\Box}\mbox{}\\

In order to prove \ref{propU} we need some preparations. Since $U^l_{pq}(H)$ is defined over $\overline{\mathbb Q}$ 
we may assume that $K=\mathbb C$. Let $R:=\mathbb C[\{A_{\nu\mu}\}]$ and $F\in R[X,Y]$ the generic Weierstra\ss\ 
polynomial of type $p,q$. We have $\Spec R =\mathbb A^n(\mathbb C)$ with $n=\nfrac{1}{2}(p+1)(q+1)-1$. 
In $\Spec{R[X,Y]}=\mathbb A^n(\mathbb C)\times \mathbb A^2(\mathbb C)$ we consider the smooth 
subschemes $V(F,F_X,F_Y)\cong \mathbb A^{n-1}(\mathbb C)$ and $V(F_X,F_Y)\cong \mathbb A^n(\mathbb C)$. 
Let $R'=R/(\Delta)$ be the image of $R$ in $R[X,Y]/(F,F_X,F_Y)$ and
$$
	\mathcal H^l_{pq}\subset \mathcal H_{pq}\subset \mathcal H=\Spec{R'}\subset \Spec R =\mathbb A^n(\mathbb C)
$$
as in Section 1. Further let $\pi: \mathbb A^n(\mathbb C)\times \mathbb A^2(\mathbb C)\to \mathbb A^n(\mathbb C)$ 
be the projection onto the first factor. Its restriction $\pi_0: V(F,F_X,F_Y)\to \mathbb A^n(\mathbb C)$ to 
$V(F,F_X,F_Y)$ is finite and has image $\mathcal H$. For a closed point $\alpha\in \mathcal H^l_{pq}$ 
the corresponding curve $C_{\alpha}$ has $l$ nodes $(x_1,y_1),\dots,(x_l,y_l)$ and no other singularities.

We endow $\mathbb C^m\ (m>0)$ with its standard norm \norm{ } and standard topology. 
For $P\in \mathbb C^m$ and $\epsilon>0$ let 
$U_{\epsilon}(P):=\{Q\in \mathbb C^m\vert\ \norm{Q-P}<\epsilon\}$  
denote the $\epsilon$-neighborhood of $P$.

The proof of the following proposition is inspired by arguments of Benedetti-Risler [BR], Lemma 5.5.9 and 
Pecker [P] in real algebraic geometry.

\prop[Simplification of nodal curves] {
Let $P_{i_1},\dots,P_{i_{\lambda}}$ be distinct nodes of $C_{\alpha}\ (1\le \lambda\le l)$. 
Given $\epsilon>0$ and $\delta>0$ there exists $\beta\in U_{\epsilon}(\alpha)$ such that the curve 
$C_{\beta}: F(\beta,X,Y)=0$ has $\lambda$ distinct nodes $Q_1,\dots,Q_{\lambda}$ and no other 
singularities where $Q_k\in U_{\delta}(P_{i_k})$ for $k=1,\dots,\lambda$.
}

{
We obtain Proposition \ref{propU} by applying 3.3 to the normed curve $C_{\alpha}$ associated to the Lissajous curve 
of type $p,q$. Let $(x_i,y_i)\ (i=1,\dots,d)$ be the nodes of $C_{\alpha}$ and 
$\gamma_i=(p-1)(q-1)-1-(a_ip+b_iq)\ (i=1,\dots,d)$ the gaps of $H_{pq}$. Then the determinant
$$
	D_{\mathbb N}(x_1,y_1,\dots,x_d,y_d)=\det{\left(x_i^{a_j}y_i^{b_j}\right)_{i,j=1,\dots,d}}
$$
does not vanish by Lemma 2.1. Let $\gamma_{j_k}\ (k=1,\dots,l)$ be the gaps of $H_{pq}$ which are closed in $H$.
Consider the columns of $\left(x_i^{a_j}y_i^{b_j}\right)$ corresponding to the $(a_{j_k},b_{j_k})\ (k=1,\dots,l)$. 
Since they are linearly independent there exist nodes $P_{i_k}:=(x_{i_k},y_{i_k})$ of the curve $C_{\alpha}$ such that
$$
	D_H(x_{i_1},y_{i_1},\dots,x_{i_l},y_{i_l})\ne 0.
$$
By Proposition 3.3 there is a nodal curve $C_{\beta}:F(\beta,X,Y)=0$ with exactly $l$ nodes 
$Q_k=(\xi_k,\eta_k)\ (k=1,\dots,l)$ which are arbitrarily close to the $P_{i_k}$. Then for a suitable $\beta$ also \linebreak 
$D_H(\xi_1,\eta_1,\dots,\xi_l,\eta_l)\ne 0$, and it follows that $(\beta,\xi_1,\eta_1,\dots,\xi_l,\eta_l)$
$\in U^l_{pq}(H)$.
}\hfill \ensuremath{\Box}

\prf[Proposition 3.3]{
In the following we consider $S:=V(F,F_X,F_Y)\cap\mathbb C^n\times \mathbb C^2$ and 
$T:=V(F_X,F_Y)\cap \mathbb C^n\times\mathbb C^2$ as submanifolds of $\mathbb C^n\times \mathbb C^2$. Then 
$S\cong \mathbb C^{n-1}$ is a hypersurface in $T\cong \mathbb C^n$. We shall study the holomorphic maps 
$\pi:\mathbb C^n\times \mathbb C^2\to \mathbb C^n$ and $\pi_0: S\to \mathbb C^n$ corresponding to the 
morphisms $\pi$ and $\pi_0$ from above in the neighborhood of $\alpha\in \mathbb C^n$. We have
$$
	\pi_0^{-1}(\alpha)=\{\alpha\}\times\text{Sing}(C_{\alpha})=\{(\alpha,x_i,y_i)\vert i=1,\dots,l\}.
$$

\lma{
Given $\delta >0$ there are for small $\epsilon >0$ open neighborhoods $U_i$ of $(\alpha,x_i,y_i)$ on 
$S\ (i=1,\dots,l)$ with the following properties:
\begin{enumerate}
\item[(i)] The $U_i$ are pairwise disjoint and
	$$
		\pi_0^{-1}(U_{\epsilon}(\alpha))=\bigcup_{i=1}^lU_i,\ U_i\subset U_{\epsilon}(\alpha)\times 
		U_{\delta}(x_i,y_i)\ \text{for}\ i=1,\dots,l.
	$$
\item[(ii)] $\pi(U_i)\subset U_{\epsilon}(\alpha)$ is a submanifold of codimension $1\ (i=1,\dots,l)$ and the map 
	$\pi_0: U_i\to \pi(U_i)$ is biholomorphic.
\item[(iii)] For any subset $\{j_1,\dots,j_{\lambda}\}\subset \{1,\dots,l\}$ with $\lambda$ distinct elements 
	$\pi(U_{j_1})\cap\dots\cap\pi(U_{j_{\lambda}})$ is a submanifold of $U_{\epsilon}(\alpha)$ of codimension $\lambda$.
\end{enumerate}
}

Using the lemma we can finish the proof of Proposition 3.3 as follows: Since $\mathcal H_{pq}$ 
is open in $\mathcal H$ we can choose in 3.4 an $\epsilon>0$ such that 
$U_{\epsilon}(\alpha)\cap\mathcal H\subset\mathcal H_{pq}$. 
Then for all $\beta\in U_{\epsilon}(\alpha)\cap\mathcal H$ it follows that $C_{\beta}$ is 
a nodal curve of type $p,q$. By dimension reasons the set
$$
	B:=\pi(U_{i_1})\cap\dots\cap\pi(U_{i_{\lambda}})\setminus\bigcup_{i\not\in\{i_1,\dots,i_{\lambda}\}}\pi(U_i)
$$
is not empty. Moreover since the $U_i\subset U_{\epsilon}(\alpha)\times U_{\delta}(x_i,y_i)$ are 
pairwise disjoint and $\pi_0: U_i\to \pi(U_i)$ is bijective, for any $\beta\in B$ the fiber $\pi_0^{-1}(\beta)$ 
consists of exactly $\lambda$ points $(\beta,Q_k)\in U_{i_k}$ and $Q_k\in U_{\delta}(P_{i_k})$ for 
$k=1,\dots, \lambda$.
}

\prf[Lemma 3.4]{
(i) We shall apply the Implicit Function Theorem to the map 
$(F_X,F_Y): \mathbb C^n\times\mathbb C^2\to \mathbb C^2$ given by $F_X$ and $F_Y$. 
Remember that $S\cong \mathbb C^{n-1}$ is a hypersurface of $T=\{(\beta,x,y)\vert F_X(\beta,x,y)=F_Y(\beta,x,y)=0\}$.
The Jacobian of the map $(F_X,F_Y)$ has rank 2 at the points $(\alpha,x_i,y_i)$ since the 
Hessian $\text{Hess}_F$ is one of its 2-minors and $\text{Hess}_F(\alpha,x_i,y_i)\ne 0$ for $i=1,\dots,l$.

The Implicit Function Theorem states that there exist $\epsilon_0>0$ and $\delta_0>0$ and 
holomorphic maps $\varphi_i:U_{\epsilon_0}(\alpha)\to U_{\delta_0}(x_i,y_i)$  with 
$\varphi_i(\alpha)=(x_i,y_i)$ such that $T\cap U_{\epsilon_0}(\alpha)\times U_{\delta_0}(x_i,y_i)$ 
is the graph $\Gamma_{\varphi_i}=\{(\beta,\varphi_i(\beta))\vert \beta\in U_{\epsilon_0}(\alpha)\}$ 
of $\varphi_i\ (i=1,\dots,l)$.
The morphism $\pi_0$ of $\mathbb C$-schemes is finite. Then the underlying continuous map $\pi_0$ 
is closed with respect to the standard topology, as is well-known.
Further $U:=\bigcup_{i=1}^lS\cap \Gamma_{\varphi_i}$ is an open neighborhood of $\pi_0^{-1}(\alpha)$ on $S$.
Hence $W:=\mathbb C^n\setminus\pi_0(S\setminus U)$ is an open neighborhood of $\alpha$ in 
$\mathbb C^n$ such that $\pi_0^{-1}(W)\subset U$. For small $\epsilon\le \epsilon_0$ we 
have $\pi_0^{-1}(U_{\epsilon}(\alpha))\subset U$ and so 
$$
	\pi_0^{-1}(U_{\epsilon}(\alpha))=\bigcup_{i=1}^lU_i
$$
where
$U_i:=S\cap(\Gamma_{\varphi_i}\cap \pi^{-1}(U_{\epsilon}(\alpha)))
=S\cap\Gamma_{\varphi_i\vert U_{\epsilon}(\alpha})$
is an open neighborhood of $(\alpha,x_i,y_i)$ on $S\ (i=1,\dots,l)$.
For small $\epsilon>0$, as $\varphi_1,\dots,\varphi_l$ are continuous functions, 
the $U_1,\dots,U_l$ are pairwise disjoint and 
$U_i\subset U_{\epsilon}(\alpha))\times U_{\delta}(x_i,y_i)$ for $i=1,\dots,l$.

\noindent
(ii) Since $U_i\subset \Gamma_{\varphi_i\vert U_{\epsilon}(\alpha)}$ is a submanifold 
of codimension $1$ and $\pi: \Gamma_{\varphi_i\vert U_{\epsilon}(\alpha)}\to U_{\epsilon}(\alpha)$ 
is biholomorphic $\pi(U_i)\subset U_{\epsilon}(\alpha)$ is likewise a submanifold of 
codimension $1$ and $\pi:U_i\to \pi(U_i)$ is biholomorphic.

\noindent
(iii) The gradient of $F$ at $(\alpha,x_i,y_i)$ has the form $(v_i,0,0)$ with 
$v_i:=(\{x_i^{\nu}y_i^{\mu}\}_{\nu p+\mu q<pq})$ for $i=1,\dots,l$.
By 2.1 the vectors $v_i$ are linearly independent, and $v_i$ is normal to the hypersurface 
$\pi(U_i)$ at $\alpha$. It  follows that $\pi(U_{i_1})\cap\dots\cap\pi(U_{i_{\lambda}})$ 
is for small $\epsilon>0$ a submanifold of $U_{\epsilon}(\alpha)$ of codimension $\lambda$.
}\mbox{}\\

In connection with Theorem \ref{thmFamW} we have a question: Given a Weierstra\ss\ semigroup close its greatest gap. 
Do we get again a Weierstra\ss\ semigroup?
\vskip3mm
\noindent
{\bf Acknowledgement.} We are grateful to Reinhold H\"ubl for valuable discussions and hints.

\end{document}